%
\input amstex
\documentstyle {amsppt}
\loadbold
\font\smbf=cmbx10 at 11pt
\font\medbf=cmbx10 scaled \magstep1
\font\bigbf=cmbx10 scaled \magstep2
\magnification=1000
\TagsOnRight
\input psfig.sty
\hsize=6in
\vsize=8in
\hfuzz=3 pt
\hoffset=0.1in
\tolerance 9000
\def\L{\text{{\bf L}}}
\def\meas{\hbox{meas}}
\def\vvs{\vskip 0.5em}
\def\v{\vskip 1em}
\def\vs{\vskip 2em}
\def\vsk{\vskip 4em}
\def\a{\alpha}
\def\tv{\hbox{\text{Tot.Var.}}}

\def\bn{\text{\bf n}}
\def\la{\big\langle}

\def\ra{\big\rangle}
\def\ov{\overline}

\def\R{\Bbb R}
\def\A{{\Cal A}}

\def\K{{\Cal K}}
\def\BV{{\strut \text{{\rm BV}}}}
\def\ds{\diamondsuit}

\def\Om{\Omega}

\def\forallt{\hbox{for all~}}

\def\ve{\varepsilon}
\def\eps{\varepsilon}
\def\n{\noindent}
\def\La{\Big\langle}
\def\Ra{\Big\rangle}
\def\fine{\hfill $\square$}

%
\voffset=0.4in
%
%
\def\AB{\hbox{[A-B]}}
%

%

%
%
%
\def\Bimp{\hbox{[B1]}}
%
%
\def\Bro{\hbox{[Bro]}}
%
\def\CLSS{\hbox{[C-L-S-S]}}
%

%
%

%
%
\def\CLRSfslf{\hbox{[C-L-R-S]}}
%

%
%
\def\Cornc{\hbox{[Cor1]}}
%

%
\def\Cortv{\hbox{[Cor3]}}
%

%

%

%
%
\def\Hdvf{\hbox{[He1]}}
%
%

%
%

%
%
\def\KS{\hbox{[K-S]}}
%
%
\def\LSrrs{\hbox{[L-S1]}}
%
%
\def\LScrs{\hbox{[L-S2]}}
%
%

%
%
%
\def\Riscl{\hbox{[Ri1]}}
%
%
%
\def\Risclsf{\hbox{[Ri2]}}
%
\def\Ry{\hbox{[Ry]}}
%
%

%

%

%
\def\Ss{\hbox{[So4]}}
%
\def\SS{\hbox{[SS]}}
%
\def\Su{\hbox{[Su]}}

\def\rhead=\vbox to 2truecm{\hbox
{\ifodd\pageno\rightheadline\else\leftheadline\fi}
\hbox to 1in{\hfil}\psill}
\tenrm
\def\rightheadline{\hfil{\eightrm Flow Stability
of Patchy Vector Fields
}\hfil\folio}
\def\leftheadline{\folio\hfil{\eightrm F. Ancona and A. Bressan}\hfil}
\def\onepageout#1{\shipout\vbox
\offinterlineskip
\vbox to 2in{\rhead}
\vbox to \pageheight
\advancepageno}
\NoBlackBoxes
%
\document
\topskip=20pt
\nopagenumbers

\null
\v
\vskip 2cm
\centerline{\bigbf Flow Stability of Patchy Vector Fields}
\v
\centerline{\bigbf and Robust Feedback Stabilization}
\vskip 2cm 
\centerline{\it Fabio Ancona$^{(*)}$ and Alberto Bressan$^{(**)}$}
\vs

\parindent 40pt

\item{(*)} Dipartimento di Matematica and C.I.R.A.M., Universit\`a di Bologna,
\item{} Piazza Porta S.~Donato~5,~Bologna~40127,~Italy.
\item{} e-mail: ancona\@ciram3.ing.unibo.it.
\v
\item{(**)} S.I.S.S.A., Via Beirut~4,~Trieste~34014,~Italy.
\item{} e-mail: bressan\@sissa.it.
\vskip 1cm

\centerline{June 2001}
\vskip 2cm

\n {\bf Abstract.} 
The paper is concerned with {\it patchy vector fields}, a class
of discontinuous, piecewise smooth vector fields
that were introduced in \AB \
to study feedback stabilization problems.
We prove the stability of the corresponding
solution set w.r.t.~a wide class
of impulsive perturbations.
These results yield the robusteness of
{\it patchy feedback controls}
in the presence of measurement errors
and external disturbances.
\vs
\vs

\noindent
{\it 1991 AMS-Subject Classification.} \
%
%
34 A,
%
%
34 D, 
%
%
49 E, 
%
%
93 D.
\vfill\eject
\vsk

\pageno=1

\parindent 20pt

\n{\medbf 1 - Introduction and Basic Notations.}
\v
Aim of this paper is to establish the stability of the set of 
trajectories
of a patchy vector field w.r.t.~various types of perturbations, 
and the robustness of patchy feedback controls.

Patchy vector fields were introduced in \AB\ in order to study 
feedback stabilization problems.  
The underlying motivation is the following: The analysis of 
stabilization problems by means of Lyapunov functions 
usually leads to stabilizing feedbacks with a wild set of
discontinuities.  
On the other hand, as shown in \AB, by patching together
open loop controls one can always construct a piecewise smooth
stabilizing feedback whose discontinuities have a very simple structure.
In particular, one can develop the whole theory by studying
the corresponding discontinuous O.D.E's within the
classical framework of Carath\'eodory solutions. 
We recall here the main definitions:
\v
\n{\bf Definition 1.1.} By a {\it patch} we mean a pair 
$\big(\Omega,\, g\big)$ where 
$\Omega\subset\R^n$ 
is an open domain with smooth boundary $\partial\Om,$ and $g$ is a 
smooth vector field 
defined on a neighborhood of the closure
$\ov\Omega$, which points strictly inward at each
boundary point $x\in\partial\Om$.
\v  
Calling $\bn(x)$ the outer normal at the boundary
point $x$, we thus require 
$$
\la g(x),~\bn(x)\ra <0\qquad\forallt\quad x\in\partial\Om.
\tag 1.1
$$
\v
\n{\bf Definition 1.2.}  \ We say that $g:\Omega\mapsto\R^n$
is a {\it patchy vector field} on the open domain $\Om$
if there exists a family of patches
$\big\{ (\Omega_\alpha,~g_\alpha) ;~~ \alpha\in\A\big\}$ such that

\n - $\A$ is a totally ordered set of indices;

\n - the open sets $\Omega_\alpha$ form a locally 
finite covering of $\Omega$, i.e. $\Omega=\cup_{\alpha\in \A}\Omega_\alpha$
and every compact set~$K\subset \Bbb R^n$
intersect only a finite number of domains $\Omega_\alpha, \, \alpha \in \A$;

\n - the vector field $g$ can be written in the form
$$
g(x) = g_\alpha (x)\qquad \hbox{if}\qquad x \in 
\Omega_\alpha \setminus 
\displaystyle{\bigcup_{\beta > \alpha} \Omega_\beta}.
\tag 1.2               
$$
\v
By setting
$$
\alpha^*(x) \doteq \max\big\{\alpha \in \A~;~~ x \in \Omega_\alpha
\big\},
\tag 1.3
$$
we can write (1.2) in the equivalent form
$$
g(x) = g_{_{\alpha^*(x)}}(x) \qquad \forallt~~x \in \Omega.
\tag 1.4
$$
We shall occasionally adopt the longer notation
$\big(\Omega,\ g,\ (\Omega_\alpha,\,g_\alpha)_{_{\alpha\in \A}} \big)$ 
to indicate a patchy vector field, specifying both the domain
and the single patches.
If $g$ is a patchy vector field, the differential equation
$$
\dot x = g(x)
\tag 1.5
$$
has many interesting properties.
In particular, in \AB\ it was proved that the set 
of Carath\'eodory solutions of
(1.5) is closed in the topology of uniform convergence, but possibly not
connected.  
Moreover, given an initial condition
$$
x(t_0)=x_0,
\tag 1.6
$$
the Cauchy problem (1.5)-(1.6) has at least one forward
solution, and at most one backward solution.
For every  Carath\'eodory solution \, $x=x(t)$ \, of (1.5),
the map \, $t \mapsto \a^*(x(t))$ \, is left continuous
and non-decreasing.

In this paper we study the stability of the solution set for (1.5) 
w.r.t.~various perturbations. 
Most of our analysis will be concerned with
impulsive perturbations, described by
$$
\dot y = g(y)+ \dot w.
\tag 1.7
$$
Here $w=w(t)$ is any 
left continuous function with bounded variation.
By a solution of the perturbed system (1.7) 
with an initial condition
$$
y(t_0)=y_0,
\tag 1.8
$$
we mean
a measurable function $t\mapsto y(t)$ such that
$$
y(t)=y_0+\int_{t_0}^t g\big(y(s)\big)\,ds +
\big[w(t)-w(t_0)\big].
\tag 1.9
$$
If $w(\cdot)$ is discontinuous, the system (1.7) has impulsive
behavior and the solution $y(\cdot)$ will be discontinuous as well.
We choose to work with (1.7) because it provides a simple and general
framework to study robustness properties.  Indeed, consider a system with
both inner and outer perturbations, of the form
$$
\dot x=g\big( x + e_1(t)\big) + e_2(t).
\tag 1.10
$$
The map $t\mapsto 
y(t)\doteq x(t)+e_1(t)$ then satisfies the impulsive equation
$$
\dot y=g(y)+e_2(t)+\dot e_1(t)=g(y)+\dot  w,
$$
where
$$ 
w (t)=e_1(t)+\int_{t_0}^t e_2(s)\,ds.
$$
Therefore, from the 
stability of solutions of (1.7) w.~r.~t.~small BV
perturbations $w$, one can immediately \, deduce a 
result on the \, stability of solutions of (1.10), when 
\, $\tv\{e_1\}$ \,
and \, $\|e_2\|_{\L^1}$ \, are suitably 
small.
Here, 
\, $\tv\{e_1\}$ \, denotes the total variation
of the function \, $e_1$ \, over the whole interval where it is defined,
while $\tv \big\{e_1\,;\,J\,\big\}$ \, 
denotes
the total variation of  \, $e_1$ \,
over a  subset~$J$.\, 
Any BV function $w=w(t)$ can be redefined up to 
$\L^1$-equivalence.
For sake of definiteness, throughout the paper
we shall always consider
left continuous representatives,
so that 
$w(t)=w(t^-)\doteq \!\!\!\underset{\ \ s\to t-}\to \lim w(s)$
for every $t$.
The Lebesgue measure of a Borel set $J\subset \Bbb R$
will be denoted by $\text{\rm meas}(J)$.
\v
We observe that,
since the Cauchy problem for (1.5) does not have 
forward uniqueness and continuous dependence, one clearly cannot
expect that a single solution of (1.5) be stable under small
perturbations.  
What we prove is a different stability property, involving
not a single trajectory but the
whole solution set: If the perturbation $w$ is
small in the BV norm, then every solution of (1.7) is close to some
solution of (1.5). This is essentially an upper semicontinuity 
property of the solution set.
\vs
\n{\bf Theorem 1.} {\it Let $g$ be a 
patchy
vector field on an open domain $\Omega\subset \R^n.$
Consider a sequence of solutions $y_\nu(\cdot)$
of the perturbed system
$$
\dot y_\nu = g(y_\nu)+ \dot 
w_\nu\qquad\qquad t\in[0,T],
\tag 1.11
$$
with \, $\tv\{w_\nu\}\to 0$ \, as \, $\nu\to\infty$. \, 
If the $y_\nu : [0,T]\mapsto\Omega$ converge
to a function $y: [0,T]\mapsto\Omega,$ uniformly
on $[0,T]$, 
then $y(\cdot)$ is a Carath\'eodory solution of (1.5)
on $[0,\,T].$}
\vs

\n{\bf Corollary 1.3.} {\it Let $g$ be a 
patchy vector field
on an open domain $\Omega\subset \R^n.$
Given any closed subset $A\subset \Omega,$ 
any compact set $K\subset A,$ and any $T, \ve>0$, 
there exists \,
$\delta=\delta(A,K,T,\ve)>0$ \, such that the following holds.
If $y:[0,T]\mapsto A$ is a solution of the perturbed system (1.7),
with $y(0)\in K$ and $\tv \{w\}<\delta$, then there exists a solution
$x:[0,T]\mapsto \Omega$ of the unperturbed equation (1.5) with}
$$
\big\|x-y\big\|_{\L^\infty([0,T])}
<\ve\,.\tag 1.12
$$
\v

We remark that the type of stability described above is precisely what
is needed in many applications to feedback control.
As an example, consider the problem of stabilizing to the origin
the control system
$$
\dot x=f(x,u).\tag 1.13
$$
Given a compact set $K$ and $\ve>0$,
assume that there exists a piecewise constant feedback $u=U(x)$ such that
$g(x)\doteq f\big(x,U(x)\big)$ is a patchy vector field, and such that
every solution of (1.5) starting from a point $x(0)\in K$ 
is steered inside the ball $B_\ve$ centered at the origin with radius
$\ve$, within a time $T>0$.
By Corollary 1, if the perturbation $w$ is sufficiently small (in the
BV norm), every solution of the perturbed system (1.7) will be
steered inside the ball $B_{2\ve}$ within time $T$.
In other words, the feedback still performs well in the presence of
small perturbations.
Applications to feedback control will be discussed in more detail in
Section~3.


\vsk
\n{\medbf 2 - Stability of Patchy Vector Fields.}
\v
We begin by proving a
local existence result for solutions
of the perturbed system (1.7).
\v
\n{\bf Proposition 2.1.} {\it Let $g$ be a 
patchy
vector field on an open domain $\Omega\subset \R^n.$
Given any compact set $K \subset \Omega$,
there exists  $\overline\chi=\overline\chi_K>0$ such that,
for each  $y_0\in K, \, t_0\in\Bbb R,$ \,
and for every  Lipschitz continuous function
$w=w(t),$ with Lipschitz constant $\|\dot w\|_{\L^\infty}<\overline\chi$,
the Cauchy problem (1.7)-(1.8) has at least one local
forward solution.}
\v
\noindent
{\bf Proof.} \ 

%
%

\noindent
Fix some compact subset $K'\subset\Omega$
whose interior contains $K$.
To prove the
local existence  of a forward solution
to (1.7), first observe that,
because of the inward-pointing condition
condition (1.1) and the smoothness assumptions
on the vector fields $g_\alpha$, one can find for any $\alpha\in \Cal A$
some constant $\chi_\alpha>0$ such that
$$
\sup
\Sb
\\
\noalign{\vskip -4pt}
\\
x \in \partial \Omega_\alpha\cap K'
\\
|v|\leq \chi_\alpha
\endSb
\big\langle g_\alpha(x) + v,\,\bn_{\boldsymbol\alpha}(x)\big\rangle < 0,
\tag 2.1
$$
where
$\bn_{\boldsymbol\alpha}(x)$ is the outer 
normal to $\partial \Omega_\alpha$
at the boundary
point $x$. Since $K'$ is a compact set
and $\{\Omega_\alpha\}_\alpha$ is
a locally finite covering of $\Omega$, there will be 
only finitely many elements of $\{\Omega_\alpha\}_\alpha$
that intersect $K'.$  Let
$$
\big\{\,
{\alpha_1},\dots,{\alpha_{\!_N}}
\,\big\}
=
\big\{\,
\alpha\in \Cal A~:~\Omega_\alpha \cap K' \neq \emptyset
\,\big\},
\tag 2.2
$$
and, by possibly renaming the indices $\alpha_i,$
assume  that
$$
\alpha_1< \dots < \alpha_{\!_N}.
\tag 2.3
$$
Choose a constant $\overline\chi>0$
such that
$$
\overline\chi \leq \inf \big\{\chi_{\alpha_i}~:~ 
i=1,\dots,N\big\}.
\tag 2.4
$$
For any 
fixed $y_0\in K$, consider the index
$$
\widehat\a(y_0)\doteq\max\big\{ \alpha~:~y_0\in\overline\Omega_\a\big\}.
$$
By the definition of $\overline\chi,$ any solution $y=y(\cdot)$
to the Cauchy problem 
$$
\dot y = g_{\widehat\a}(y) + \dot w,
\qquad\quad y(t_0)=y_0,
$$
associated to a piecewise Lipschitz map $w=w(t)$ with 
$\|\dot w\|_{\L^\infty}<\overline\chi,$ 
remains inside $\Omega_{\widehat\a}$ for all
$t\in [t_0, t_0+\delta],$ for some $\delta>0.$
Hence, it provides also a solution to (1.6)
on some interval $[t_0, t_0+\delta'],$ \, $0<\delta'\leq \delta.$
\fine
\vs
Toward a proof of Theorem 1, we first derive
an intermediate result.
By the basic properties of a patchy vector field,
for every solution $t\mapsto x(t)$ of (1.5) the corresponding map
$t\mapsto \alpha^*\big(x(t)\big)$ in (1.3) is nondecreasing.
Roughly speaking, a trajectory can move from
a patch $\Omega_\alpha$ to another patch $\Omega_\beta$ only if
$\alpha<\beta$. This property no longer holds in the presence of
an impulsive perturbation. However, the next proposition shows
that for a solution $y$ of (1.7) the corresponding map
$t\mapsto \alpha^*\big(y(t)\big)$ is still nondecreasing, 
after a possible modification on a
small set of times.  Alternatively, one can slightly modify
the impulsive perturbation $w$, say replacing it by another
perturbation $w^\ds$, such that the map 
$t\mapsto\alpha^*\big(y^\ds(t)\big)$ is monotone
along the corresponding trajectory $t\mapsto y^\ds(t)$.
\v
\n{\bf Proposition 2.2.} {\it Let $g$ be a 
patchy
vector field on an open domain $\Omega\subset \R^n,$
determined by the family of patches
$\big\{ (\Omega_\alpha,~g_\alpha)  ; \ \alpha\in\Cal A\big\}$.
For any $T>0$ and any compact set $K \subset \Omega$, 
there exist
constants \, $C,\delta>0$ and an integer $N$
such that the following holds.

\item{(i)} \  For every $w\in\BV$ with \, $\tv\{w\}<\delta,$ \,
and for every solution  \, $y:[0,T]\mapsto \Omega$ \, of the
Cauchy problem (1.7)-(1.8) with \, $y_0\in K,$ \,
there is a  partition of \, $[0,\,T],$ \
$0=\tau_1\leq \tau_2\leq \dots \leq\tau_{N+1}=T,$ 
and indices
$$
\alpha_1<\alpha_2<\cdots <\alpha_N,
\tag 2.5
$$
%
such that
$$
\alpha^*(y(t))\geq \alpha_i
\qquad\quad \forall~t\in ]\tau_i,\,\tau_{i+1}],
\qquad\quad
i=0,\dots,N,
\tag 2.6
$$
$$
\text{\rm meas}\,
\Big(\bigcup_{i\geq 0}
\big\{
t\in [\tau_i,\,\tau_{i+1}]~:~ \alpha^*(y(t)) > \alpha_i
\big\}\Big)
< C \cdot \text{\rm Tot.Var.}\{w\}.
\tag 2.7
$$
%
%

\v
\item{(ii)} \ 
For every BV function \, $w=w(t)$ \, with \, $\tv\{w\}<\delta,$ \,
and for every solution \, $y:[0,T]\mapsto \Omega$  \, of the
Cauchy problem (1.7)-(1.8) with \, $y_0\in K,$ \,
there is a BV function \, $w^\ds=w^\ds(t)$ \,
and a solution \,
$y^\ds:[0,T]\mapsto \Omega$ \, of 
$$
\dot y^\ds=g(y^\ds)+\dot w^\ds\,,
\tag 2.8
$$
so that the map $t \mapsto \alpha^*(y^\ds(t))$ is 
non-decreasing and left continuous, and there holds
$$
\aligned
\tv\{w^\ds\}
&\leq C\cdot \tv\{w\}\,,
\\
\noalign{\medskip}
\big\|y^\ds-y\big\|_{\L^\infty([0,T])}
&\leq C\cdot \tv\{w\}\,.
\endaligned
\tag 2.9
$$
}
\v
\noindent
{\bf Proof.} \

\noindent
{\smbf 1.} 
The proof of (i) will be given in three steps.
\v
\noindent
{\smc Step 1.} \ Since each $g_\alpha$ is a smooth
vector
field and we are assuming a uniform bound on the total variation
of every perturbation $w=w(t)$, there will be some compact
subset $K'\subset\overline\Omega$ that contains
every solution $y : [0,T]\mapsto\Omega$ of (1.7)
starting at a point $y_0\in K$.
We will assume without loss of generality that every 
domain \, $\Omega_\alpha$ \, is bounded since, otherwise,
one can replace \, $\Omega_\alpha$ \, with its intersection 
\, $\Omega_\alpha \cap \Omega'$ \, with 
a bounded domain \, $\Omega'\subset \Omega$ \, that contains \, $\overline K'$,
preserving the inward-point condition (1.1).
For each $\alpha\in\Cal A$, define the map
$\varphi_\alpha:\Omega\mapsto \Bbb R$ by setting
$$
\varphi_\alpha(x)\doteq
\cases
d(x,\,\partial \Omega_\alpha)
\quad\ &\text{if}\qquad x\in \Omega_\alpha,
\\
\noalign{\smallskip}
-d(x,\,\partial \Omega_\alpha)
\quad\ &\text{otherwise,}
\endcases
\tag 2.10
$$
and let 
$$
\varphi_\alpha^+(x)\doteq\max\{\varphi_\alpha(x),\, 0\}
$$
denote the positive part of \, $\varphi_\alpha(x).$ \,
The regularity assumptions on the patch $\Omega_\alpha$
guarantee that  $\varphi_{\alpha}$
is smooth if restricted to a sufficiently small neighborhood
of the boundary $\partial \Omega_{\alpha}$. 
Thus, if 
$\{\Omega_{\alpha_i}~:~ 
i=1,\dots,N\}$ denotes the finite collection of domains
that intersect the compact set $K'$ as in~(2.2)-(2.3),
there will be some constant $\overline\rho>0$
so that, setting 
$$
\Omega_\alpha^{\overline\rho}\doteq
\big\{
x\in\Omega
~:~ d(x,\,\partial \Omega_\alpha)\geq \overline\rho
\big\},
\tag 2.11
$$
%
the restriction of any map $\varphi_{\alpha_i}$ to the domain 
$\Omega\setminus \Omega_{\alpha_i}^{\overline\rho}$
be smooth.
In particular, for any $i=1,\dots,N,$
we will have
$$
\nabla \varphi_{\alpha_i}(x)=-\bn_{{\boldsymbol\alpha}_{\bold i}}
\big(\pi_{\alpha_i}(x)\big),
\qquad\quad \forall~~~x\in
\Omega\setminus 
\Omega_{\alpha_i}^{\overline\rho},
\tag 2.12
$$
where $\bn_{{\boldsymbol\alpha}_{\bold i}}$ represents
as usual the outer normal to $\partial \Omega_{\alpha_i},$ while
$\pi_{\alpha_i}(x)$ denotes the projection of the point $x$
onto the set $\partial \Omega_{\alpha_i}.$
On the other hand, thanks to the inward-pointing condition (1.1),
we can choose the constant ${\overline\rho}$ so that
$$
\sup
\Sb
\\
\noalign{\vskip -4pt}
\\
i=1,\dots,N
\\
x \in \Omega_{\alpha_i}\setminus \Omega_{\alpha_i}^{\overline\rho}
\endSb
\big\langle g_{\alpha_i}(x),\,\bn_{{\boldsymbol\alpha}_{\bold i}}
(\pi_{\alpha_i}(x))\big\rangle \leq -c',
\tag 2.13
$$
for some $c'>0$.  Moreover, the smoothness  of the
fields $g_\alpha$ on $\overline \Omega$
implies the existence of some $c''>0$
such that
$$
\sup
\Sb
\\
\noalign{\vskip -4pt}
\\
i=1,\dots,N, \ 
j>i
\\
x \in \Omega_{\alpha_i}
\endSb
\Big|
\big\langle g_{\alpha_j}(x),\,\bn_{{\boldsymbol\alpha}_{\bold i}}
(\pi_{\alpha_i}(x))\big\rangle 
\Big|\leq c''.
\tag 2.14
$$

\v
\noindent
{\smc Step 2.} \ Consider now a left continuous BV function
$w=w(t)$ 
and let 
$y:[0,T]\mapsto \Omega$ be the solution
of the corresponding Cauchy problem (1.7)-(1.8), with $y_0\in K.$
Observe that, for any $i=1,\dots,N,$ and for any interval 
$J\subset [0,T]$ such that 
$$
y(t)\in \Omega\setminus \Omega_{\alpha_i}^{\overline\rho}
\qquad\quad \forallt \qquad t\in J,
$$
the composed map
$\varphi_{\alpha_i}^+ \circ y:J\mapsto \Bbb R$ 
is also a left continuous BV function
whose distributional derivative \,
$\mu_i\doteq D \,(\varphi_{\alpha_i}^+\circ y)$ \, is a Radon measure, 
which can be decomposed into an absolutely continuous $\mu_i^{ac}$
and a singular part $\mu_i^{s}$, 
w.~r.~t.~the Lebesgue measure \, $dt$. \,
One can easily verify that, for any Borel set $E\subset J$,
the absolutely continuous part of $\mu_i$
is given by
$$
\mu_i^{ac}(E)=
\int_{E^+} \big\langle \nabla \varphi_{\alpha_i}(y(t)),\ 
g(y(t))
+\dot w(t)
\big\rangle~dt, 
\qquad\quad E^+\doteq \{t\in E~;~ y(t)\in \Omega_{\alpha_i}\}\,.
\tag 2.15
$$
Moreover, calling \, $\mu_w^{ac},\,\mu_w^s,$ \, 
respectively the
absolutely continuous and the singular part of \,
$\mu_w\doteq \dot w,$ \,
the following bounds hold
%
$$
\align
\bigg|\int_{E^+} \big\langle \nabla \varphi_{\alpha_i}(y(t)),\ 
\dot w(t)
\big\rangle~dt 
+\mu_i^{s}(E)\bigg|&\leq
c'''\cdot
\Big\{\big|\mu_w^{ac}(E)\big|+\big|\mu_w^{s}(E)\big|
\Big\}
\\
\noalign{\smallskip}
&\leq c'''\cdot
\tv\{w\},
\tag 2.16
\endalign
$$
for some constant $c'''>0$ 
that depends only on the 
compact set $K'$ and on the time interval $[0,T]$.
Let \, $C_i, \, \ell_i, \ i=N, N-1, \dots ,1$, \, be the 
constants recursively defined by
$$
\alignat 4
C_N
&\doteq 1+c''',\qquad\qquad
&\ell_N
&\doteq \frac{2 C_N}{c'},
\tag 2.17
\\
\noalign{\smallskip}
C_i
&\doteq
c''\cdot \ell_{i+1}+\sum_{j=i+1}^N  C_j,
\qquad\qquad
&\ell_i
&\doteq 
\frac{1}{c'}\Big(2 C_i+c''\cdot \sum_{j=i+1}^N \ell_j\Big)
\qquad\text{if}\qquad i<N.
\tag 2.18
\endalignat
$$
\v

\noindent
{\bf Lemma 2.3.} {\it Assume that
$$
\tv\{w\}<
\delta\doteq
\frac{\overline\rho}{2 C_1},
\tag 2.19
$$
and assume that there exists some interval $[t_1,\,t_2]\subset [0,\,T]$
and some index $i\in\{1,\dots,N\}$ such that
$$
\meas
\big\{t\in[t_1,\,t_2]~:~
\alpha^*(y(t))=\alpha_j
\big\}\leq \ell_j\cdot \tv\{w\}
\qquad\quad\forall~j>i
\eqno(2.20)_i
$$
together with one of the following two conditions 
\v

\parindent 40pt 
\item{${{\boldkey (}\bold a_{\bold i}{\boldkey )}}$}
$$
\varphi_{\alpha_i}(y(t))<2 C_i\cdot \tv\{w\}
\qquad\forall~t\in[t_1,\,t_2],
\eqno(2.21)_i
$$
%

%
$$
\meas
\big\{t\in[t_1,\,t_2]~:~
\alpha^*(y(t))=\alpha_i
\big\}>\ell_i\cdot \tv\{w\}.
\eqno(2.22)_i
$$
\vs

\item{${{\boldkey (}\bold b_{\bold i}{\boldkey )}}$}
There exists $\tau\in[t_1,\,t_2]$ such that
$$
\varphi_{\alpha_i}(y(\tau))\geq 2 C_i\cdot \tv\{w\}.
\eqno(2.23)_i
$$
\v
\vs
\parindent 20pt

\noindent 
Then one has
$$
\varphi_{\alpha_i}(y(t_2))\geq C_i\cdot \tv\{w\}.
\eqno(2.24)_i
$$
}
\v

Towards a proof of the lemma, observe first that
the recursive definition (2.17)-(2.18) of the constants \ $C_i,\,\ell_i,$ \
and the bound (2.19) clearly imply
$$
\gather
C_i\geq  1 +c'''+c''\cdot \sum_{j=i+1}^N \ell_j,
\tag 2.25
\\
\noalign{\smallskip}
2C_i\cdot \tv \{w\}< \overline \rho.
\tag 2.26
\endgather
$$
Assume now that $(2.20)_i-(2.22)_i$ hold.
%
%
Then, using (2.13)-(2.16) and recalling (2.25)-(2.26),  we obtain
$$
\align
\varphi_{\alpha_{i}}^+(y(t_2))&\geq
\varphi_{\alpha_{i}}^+(y(t_1))
+\!\!
\int
\limits_{\{t\in[t_1,t_2]~:~\alpha^*(t)=\alpha_i\}}
\!\!\!
\big\langle \nabla \varphi_{\alpha_{i}}(y(t)),\ g_{\alpha_{i}}(y(t))
\big\rangle~dt - 
\\
&\qquad -\!\sum_{j=i+1}^N\
\int
\limits_{\{t\in[t_1,t_2]~:~\alpha^*(t)=\alpha_j\}}
\!\!\!
\Big|
\big\langle \nabla \varphi_{\alpha_{i}}(y(t)),\ g_{\alpha_j}(y(t))
\big\rangle
\Big|~dt 
-c'''\cdot\tv \{w\}
\\
&\geq
\!\!
\int
\limits_{\{t\in[t_1,t_2]~:~\alpha^*(t)=\alpha_i\}}
\!\!\!
-\big\langle 
\bn_{{\boldsymbol\alpha}_{\bold i}}
(\pi_{\alpha_i}(y(t))),\ g_{\alpha_i}(y(t))
\big\rangle~dt 
- \Big(c''\cdot \sum_{j=i+1}^N \ell_j
+c'''\Big)\cdot\tv \{w\}
\\
\noalign{\smallskip}
&\geq
\Big(\ell_i\cdot c'-c''\cdot \sum_{j=i+1}^N \ell_j
-c'''\Big)
\cdot\tv \{w\}
\\
\noalign{\smallskip}
&> C_i\cdot\tv \{w\},
\tag 2.27
\endalign
$$
proving $(2.24)_i.$ Next, assume that $(2.20)_i$ and $(2.23)_i$ hold,
%
and let
$$
\tau'\doteq
\sup\big\{t\in[t_1,\,t_2]~:~
\varphi_{\alpha_i}(y(t))>2 C_i\cdot \tv\{w\}
\big\}.
\tag 2.28
$$

\noindent
Clearly, the bound $(2.24)_i$ is satisfied if
$\tau'=t_2$ since the map $\varphi_{\alpha_i}$
is left continuous.
%
Next, consider the case $\tau'<t_2.$ By similar computations
as in (2.27), using (2.13)-(2.16) and thanks to $(2.20)_i,$ (2.25)-(2.26),
we get
$$
\align
\varphi_{\alpha_i}^+(y(t_2))&\geq
\varphi_{\alpha_i}^+(y(\tau'))
+\!\!
\int
\limits_{\{t\in[\tau',t_2]~:~\alpha^*(t)=\alpha_i\}}
\!\!\!
\big\langle \nabla \varphi_{\alpha_i}(y(t)),\ g_{\alpha_i}(y(t))
\big\rangle~dt - 
\\
&\qquad -\!\sum_{j=i+1}^N\
\int
\limits_{\{t\in[\tau',t_2]~:~\alpha^*(t)=\alpha_j\}}
\!\!\!
\Big|
\big\langle \nabla \varphi_{\alpha_{i}}(y(t)),\ g_{\alpha_j}(y(t))
\big\rangle
\Big|~dt 
-c'''\cdot\tv \{w\}
%
\\
\noalign{\smallskip}
&>\Big(2 C_i-1-c'''-c''\cdot \sum_{j=i+1}^N \ell_j
\Big)
\cdot\tv \{w\}
\\
\noalign{\smallskip}
&> C_i\cdot\tv \{w\},
\tag 2.29
\endalign
$$
thus concluding the proof of Lemma 2.3.

\v
\noindent
{\smc Step 3.} \ Assume that the bound (2.19) on the
total variation of $w=w(t)$ holds.
Set $\tau_1=0,$ \, $\tau_{\!_{N+1}}\doteq T,$ and define recursively
the points $\tau_N,\tau_{N-1},\dots,\tau_2,$
by setting
$$
\qquad
\tau_i\doteq
\inf\Big\{
t\in[0,\,\tau_{i+1}]~:~
\varphi_{\alpha_i}(y(s))\geq C_i\cdot \tv\{w\}
\qquad\forall~s\in [t,\,\tau_{i+1}]
\Big\}\,,
\qquad 1<i\leq N\,.
\qquad
\eqno (2.30)_i
$$
%
Applying Lemma 2.3 and proceeding
by backward induction on $i=N, N\!-\!1, \dots, 2,$ we will show that, 
for any $t<\tau_i,$ \, $i=2,\dots,N,$ 
one has
$$
\gathered
\meas
\big\{s\in[0,\,t]~:~
\alpha^*(y(s))=\alpha_i
\big\}\leq \ell_i\cdot \tv\{w\},
\\
\noalign{\medskip}
\varphi_{\alpha_i}(y(t))< 2 C_i\cdot \tv\{w\}.
\endgathered
\eqno (2.31)_i
$$
Indeed, if $(2.31)_N$ is not satisfied,
one of the two conditions 
${{\boldkey (}\bold a{\boldkey )}}_{\bold N}$ 
or ${{\boldkey (}\bold b{\boldkey )}}_{\bold N}$
must be true on some interval $[0,\,\overline t\,], \ \overline t<\tau_N.$
But then, by $(2.24)_N$, we have 
$$
\varphi_{\alpha_{\!_N}}(y(s))\geq C_N\cdot \tv\{w\}
\qquad\forall~s\in [\,\overline t,\,T],
$$
which contradicts the definition $(2.30)_N.$
On the other hand, if we assume that $(2.31)_j$ holds 
for \, $j=i+1,\dots,N,$ \,
but not for \, $j=i,$ \, then one of the two conditions 
${{\boldkey (}\bold a{\boldkey )}}_{\bold i}$ 
or ${{\boldkey (}\bold b{\boldkey )}}_{\bold i}$
must be true on some interval $[0,\,\overline t\,], \ \overline t<\tau_i.$
Moreover, the inductive assumptions $(2.31)_j, \ j>i,$
imply $(2.20)_i$ and hence, as above, thanks to $(2.24)_i$
we get
$$
\varphi_{\alpha_i}(y(s))\geq C_i\cdot \tv\{w\}
\qquad\forall~s\in [\overline t,\,T],
$$
reaching a contradiction with
the definition $(2.30)_i.$
\v

To conclude the proof of (i), observe that, 
thanks to $(2.31)_i,\ i =2,\dots,N,$ we have
$$
\meas
\big\{s\in[\tau_i,\,\tau_{i+1}]~:~
\alpha^*(y(s))>\alpha_i
\big\}
\leq \Big(\sum_{j>i}
\ell_j\Big)\cdot \tv\{w\}
\qquad\quad \forall~i\geq 1.
\tag 2.32
$$
Therefore, recalling the definitions of the map $\varphi_{\alpha_i}$
at (2.10), taking $\delta$ as in (2.19), and
$$
C\geq
(N+1)\cdot \sum_{j=1}^N \ell_j,
\tag 2.33
$$
from $(2.31)_i$ and (2.32) we deduce that
the partition \, $\tau_1=0\leq \tau_2\leq  \dots \leq \tau_{N+1}=T$ \,
of $[0,\,T],$ defined at $(2.30)_i,$ satisfies the properties  
(2.5)-(2.7).
\vs

\noindent
{\smbf 2.} 
Concerning (ii), let $C, \delta>0$ be the constants 
defined according to (i) and, given a BV function $w=w(t)$
with $\tv\{w\}<\delta,$
and a solution $y:[0,T]\mapsto \Omega$ of the
Cauchy problem (1.7)-(1.8) with $y_0\in K,$ consider
the partition \, $0=\tau_1\leq \tau_2\leq \dots \leq\tau_{N+1}=T,$ \,
of \, $[0,\,T],$ \
with the properties in~(i). Setting
%
$$
\tau'_i
\doteq
\inf\Big\{
t\in[\tau_i,\,\tau_{i+1}]~:~
\alpha^*(y(t))=\alpha_i
\Big\}\qquad\quad i=1,\dots,N,
$$
%
define the map 
$$
\tau(t)
\doteq
\cases
\qquad\quad\tau_i'
\qquad&\text{if}\qquad t\in ]\tau_i,\,\tau_i']
\\
\noalign{\medskip}
\sup\big\{
s\in[\tau',\,t]~:~\alpha^*(y(s))=\alpha_i\big\}
\qquad&\text{if}\qquad t\in ]\tau_i',\,\tau_{i+1}],
\endcases
\tag 2.34
$$
over any interval $]\tau_i,\,\tau_{i+1}],$ \, $i=1,\dots,N.$ \,
Notice that, in the particular case
where $\alpha^*(y(t))>\alpha_i$ for all 
$t\in]\tau_i,\,\tau_{i+1}],$ by
the above
definitions one has $\tau(t)=\tau_i'=\tau_{i+1}$
for any $t\in]\tau_i,\,\tau_{i+1}].$
Then, let  $y^\ds:[0,T]\mapsto \Omega$ \, be the map recursively
defined
by setting 
$$
\align
y^\ds(t)
&\doteq y(t)
\qquad\qquad \forall~t\in]\tau_N,\,T],
\tag 2.35
\\
\noalign{\bigskip}
y^\ds(t)
&\doteq
\cases
y^\ds(\tau_{_{i+1}}+)
\ \ &\text{if}\quad \tau_i'=\tau_{i+1},
\\
\noalign{\medskip}
y(\tau(t)+)
\ \ &\text{if}\quad \tau_i'<\tau_{i+1},
\quad \alpha^*(y(\tau(t)))>\alpha_i,
\\
\noalign{\medskip}
y(\tau(t))
\ \ &\text{if}\quad \tau_i'<\tau_{i+1},
\quad \alpha^*(y(\tau(t)))=\alpha_i,
\endcases
\qquad\  \forall~t\in]\tau_i,\,\tau_{i+1}],
\quad \ i<N,
\tag 2.36
\\
\noalign{\bigskip}
y^\ds(0)
&\doteq y^\ds(0^+)
\tag 2.37
\endalign
$$
and let $w^\ds=w^\ds(t)$ be the function defined 
as
$$
w^\ds(t)\doteq
y^\ds(t)-\int_0^t g\big(y^\ds(s)\big)\,ds
\qquad\quad\forall~t\in[0,\,T].
\tag 2.38
$$
Clearly $y^\ds, \, w^\ds$ are both  BV functions as well as
$y, \, w$. Moreover,    $y^\ds$ is a solution
of the perturbed equation (2.8).
By construction,
for every $1\leq i\leq N$ there holds
$$
\alpha^*(y^\ds(t))=
\cases
\quad \alpha_i \qquad &\text{if}\qquad \tau_i'<\tau_{i+1},
\\
\noalign{\medskip}
\alpha^*(y^\ds(\tau_{\!_{i+1}}^{\ \ \,+}))
\qquad &\text{if}\qquad \tau_i'=\tau_{i+1}
\endcases
\qquad\quad \forall~t\in]\tau_i,\,\tau_{i+1}].
\tag 2.39
$$
Hence the map $t\mapsto\alpha^*(y^\ds(t))$
is non-increasing and left-continuous.
Next, recalling (2.6)
and observing that
$$
\alpha^*(y(t))= \alpha_i
\qquad\Longrightarrow
\qquad
\aligned
\tau(t)&=t,
\\
\noalign{\medskip}
y^\ds(t)&=y(t),
\endaligned
\qquad\qquad\forall~t\in]\tau_i,\,\tau_{i+1}],
\tag 2.40
$$
defining
$$
\Cal I\doteq \bigcup_i
\big\{
t\in ]\tau_i,\,\tau_{i+1}]~:~ \alpha^*(y(t)) > \alpha_i
\big\},
\tag 2.41
$$
we have
$$
y(t)=y^\ds(t)
\qquad\quad \forall~t\in(0,\,T)\setminus\Cal I.
\tag 2.42
$$
On the other hand, by the above definitions,
calling \, $M\doteq \sup_{y\in \Omega} |g(y)|,$ \,
we derive
$$
\align
\big|\tau(t)-t\big|
&\leq \meas(\Cal I)
\qquad\quad
\forall~t\in \Cal I,
\tag 2.43
\\
\noalign{\bigskip}
\big|y^\ds(t)-y(t)\big|
&\leq \int_{\tau(t)}^t \big|g(y(s))\big|~ds+
\tv\big\{w\,;\,[0,t]\big\}
\\
\noalign{\smallskip}
&\leq M\cdot
\meas(\Cal I)+\tv\{w\}
\qquad\qquad
\forall~t\in \Cal I,
\tag 2.44
\endalign
$$
and
$$
\align
\Big|\tv\{y^\ds\}-\tv\{y\}\Big|
&\leq \tv\big\{y\,;\,\Cal I\big\}
\\
\noalign{\smallskip}
&\leq M\cdot
\meas(\Cal I)+\tv\{w\}.
\tag 2.45
\endalign
$$
Then, using (2.44)-(2.45), we obtain
$$
\align
\Big|\tv\{w^\ds\}-\tv\{w\}\Big|
&\leq \int_{\Cal I} \Big|\big|g\big(y^\ds(s)\big)\big|-
\big|g\big(y(s)\big)\big|\Big|\,ds+
\Big|\tv\{y^\ds\}-\tv\{y\}\Big|
\\
\noalign{\smallskip}
&\leq
M'\cdot \Big\{\meas(\Cal I)+\tv\{w\}\Big\},
\tag 2.46
\endalign
$$
for some constant $M'>0,$ depending only on the field $g.$
Hence, from (2.42), (2.44), (2.46), and applying~(2.7),
it follows that $y^\ds(\cdot)$ satisfies the
estimates in (2.9), for some constant $C'>0,$
which concludes the proof of~(ii).    
\v

We can now take  $\delta$ as in~(2.19) and choose $C>C'$ 
according to (2.33).
Both properties (i) and (ii) are then satisfied, completing the
proof of Proposition~2.2.
\fine

\vsk
\n{\smbf Proof of Theorem 1.}
\v
For a given sequence of solutions \, $y_\nu: [0,T]\mapsto\Omega$ \,
of the perturbed system (1.11)
with \, $\tv\{w_\nu\}\leq \delta_\nu$, 
$\,\delta_\nu \to 0$ as \, $\nu\to\infty$, \, 
assume that the $y_\nu(\cdot)$ converge
to a function $y: [0,T]\mapsto\Omega$ uniformly
on $[0,T]$, and that $y_\nu(0)$ belongs to some
compact set $K\subset \Omega$ for everyy $\nu$.
Thanks to property (ii) of Proposition~2.2, 
in connection with any pair 
$w_{\nu}(\cdot),$ \,
$y_{\nu}(\cdot),$
there will be a BV function $w_{\nu}^\ds(\cdot)$ 
and a solution $y_{\nu}^\ds(\cdot)$ of (2.8) that satisfy 
$$
\tv\{w^\ds_\nu\}
\leq C'\cdot \delta_\nu\,,
\qquad\qquad
\big\|y^\ds_\nu-y_\nu\big\|_{\L^\infty\big([0,T]\big)}
\leq C'\cdot \delta_\nu\,,
\tag 2.47
$$
for some constant $C'>0,$ independent of $\nu$. Moreover there exists
a partition 
\, $0=\tau_{_{1,\nu}}\leq \tau_{_{2,\nu}}\leq \dots 
\leq\tau_{_{N+1,\nu}}=T$ \,
of \, $[0,\,T],$ \
such that
$$
\alpha^*(y^\ds_\nu(t))=\alpha_i
\qquad\quad\forall~t\in]\tau_{_{i,\nu}},\,\tau_{_{i+1,\nu}}],
\qquad
i=1,\dots,N.
\tag 2.48
$$
%
Recalling  (1.4) and (1.9), 
because of (2.48) we have
$$
\gather
y^\ds_\nu(t) = y^\ds_\nu(0) ~+~ \sum_{\ell=1}^{i-1} 
\int_{\tau_{_{\ell,\nu}}}^{\tau_{_{\ell+\!1,\nu}}}
\!\!g_{_{\!\alpha_{_\ell}}}\!(y^\ds_\nu(s))~ds
~+\int_{\tau_{_{i,\nu}}}^t
\!\!\!g_{_{\!\alpha_i}}\!(y^\ds_\nu(s))~ds~+~
[w_\nu^\ds(t)-w_\nu^\ds(0)]
\\
\noalign{\smallskip}
\qquad\qquad
\forall~t\in[\tau_{_{i,\nu,}}\,\tau_{_{i+\!1,\nu}}]
,\qquad i=1,\dots ,N.
\tag 2.49
\endgather
$$
By possibly taking a subsequence,
we can assume that every sequence
$\big(\tau_{_{i,\nu}}\big)_{\nu\geq 1}$ converges to some limit point, say
$$
\overline \tau_i \doteq \lim_{\nu \to \infty}\tau_{_{i,\nu}}
\qquad\quad i=1,\dots ,N+1.
$$
We now observe that
$$
]\overline \tau_i,\,\overline \tau_{i+1}[~\subseteq
\bigcup_{\mu=1}^{\infty}~
\bigcap_{\nu=\mu}^{\infty}
]\tau_{_{i,\nu}},\,\tau_{_{i+\!1,\nu}}]
\qquad\quad \forall~i.
$$
Moreover, $(2.47)_2$ and the uniform convergence 
$y_\nu(\cdot)\to y(\cdot)$
yield
$$
\lim_{\nu\to\infty} 
\big\|y^\ds_\nu-y\big\|_{\L^\infty([0,T])} =0.
\tag 2.50
$$
%
{}From  $(2.47)_1$\,, (2.48)-(2.59) we now deduce
$$
\gathered
y(t) \in \overline\Omega_{\alpha_i}
\!\setminus \!\!\!\bigcup_{\beta > \alpha_i} 
\!\!\Omega_\beta,
\\
\noalign{\medskip}
\qquad
y(t) =  y(0) + \sum_{\ell=1}^{i-1} 
\int_{\overline \tau_\ell}^{\overline \tau_{\ell+\!1}}
g_{_{\!\alpha_{_\ell}}}\!(y(s))~ds
+\int_{\overline \tau_i}^{t}
\!\!g_{_{\!\alpha_i}}\!(y(s))~ds
\endgathered
\qquad\quad \forall~t\in\,]\overline \tau_i,\,\overline \tau_{i+1}],
\qquad \forall~i.
\tag 2.51
$$
In particular, on each interval 
$[\overline \tau_i,\,\overline \tau_{i+1}],$ 
the function $y(\cdot)$ is a classical
solution of 
$\dot y = g_{\alpha_i}(y)$ and satisfies 
$$
\dot y(s^-) =g_{\alpha_i}(y(s))
\qquad\quad \forall~s\in]\overline \tau_i,\,\overline \tau_{i+1}]. 
$$
Moreover observe that, 
because of the inward-pointing condition (1.1), 
the set $\big\{t\in [\overline \tau_i,\,\overline \tau_{i+1}]
: y(t)\in \partial\,\Omega_{\alpha_i}\big\}$
is nowhere dense in $[\overline \tau_i,\,\overline \tau_{i+1}].$
Thus, if $s$
is any point in $]\overline \tau_i,\,\overline \tau_{i+1}]$
such that $y(s)\in \partial\,\Omega_{\alpha_i},$
there will be some increasing sequence 
$(s_n)_n\subset]\overline \tau_i,\,\overline \tau_{i+1}[$ 
converging to $s$ and
such that $y(s_n)\in\Omega_{\alpha_i}$ for any $n$.
But this yields a contradiction with (1.1), because
$$
0
\leq
\lim_{n \to \infty}
\La 
\frac{y(s)-y(s_n)}{s-s_n}
,~\bn_{{\boldsymbol\alpha}_{\bold i}}\big(y(s)\big)
\Ra
=
\La 
\dot{y}(s-)
,~\bn\big(y(s)\big)
\Ra
=\La 
g_{\alpha_i}\big(y(s)\big)
,~\bn_{{\boldsymbol\alpha}_{\bold i}}\big(y(s)\big)
\Ra.
$$
Hence, recalling the definition (1.2), from (2.51) we conclude 
$$
\gathered
y(t) \in \Omega_{\alpha_i}
\!\setminus \!\!\!\bigcup_{\beta > \alpha_i} 
\!\!\Omega_\beta
\qquad\quad \forall~t\in]\overline \tau_i,\,\overline \tau_{i+1}],
\qquad i=1,\dots ,N,
\\
\noalign{\medskip}
y(t) = y(0) + \int_{0}^t g\big(y(s)\big)~ds
\qquad\quad \forall~t\in[0,\,T],
\endgathered
$$
proving that $y : [0,\,T] \mapsto \Omega$ 
is a Carath\'eodory solution
of (1.5) on $[0,\,T].$
\fine

\vsk
\n{\smbf Proof of Corollary 1.3.}
\v
Assuming that statement is false, we shall reach a contradiction. 
Fix any closed subset $A\subset \Omega,$ any
compact set $K\subset A,$ and assume that, for some
$T,\, \ve>0,$ there exists a sequence of solutions 
\, $y_\nu: [0,T]\mapsto A$ \,
of the perturbed system~(1.7),
with \, $y_\nu(0)\in K,$ \, $\tv\{w_\nu\}\leq \delta_\nu,$ \, 
$\delta_\nu \to 0$ as \, $\nu\to\infty,$ \, such that
the following property holds.
\v

\parindent 30pt 
\item{${{\boldkey (}\bold P{\boldkey )}}$}
\ Every solution $x: [0,T]\mapsto\Omega$ of
the unperturbed equation (1.5) satisfies
$$
\big\|x-y_\nu\big\|_{\L^\infty([0,T])}
\geq \ve\,
\qquad\quad \forall~\nu.
\tag 2.52
$$
\v

\parindent 20pt

\noindent
For each $\nu,$ call \, $y^\ds_\nu: [0,T]\mapsto\Bbb R^n$ \, the
polygonal curve with vertices at the points \,
$y_\nu(\ell \delta_\nu), \, \ell\geq 0,$
defined by setting
$$
\gather
y^\ds_\nu(t)\doteq y_\nu\big(\ell \delta_\nu\big)+
\frac{t-\ell \delta_\nu}{\delta_\nu}\cdot
\Big(y_\nu\big((\ell+1) \delta_\nu\big)-y_\nu\big(\ell \delta_\nu\big)\Big)
\\
\noalign{\medskip}
\qquad\quad\forall~t\in[\ell \delta_\nu,\,(\ell+1)\delta_\nu]\cap[0,T],
\qquad 0\leq \ell \leq \lfloor T/\delta_\nu\rfloor\,,
\tag 2.53
\endgather
$$
where $\lfloor T/\delta_\nu \rfloor$ 
denotes the integer part of \, $T/\delta_\nu.$ \,
Since every $y_\nu(\cdot)$ is a BV function that
solves the equation (1.7), it follows that
there will be some constant $C>0,$ independent on $\nu,$
such that
$$
\tv\{y_\nu~;~J\}
\leq C \cdot\meas(J)+\tv\{w_\nu~;~J\}
\tag 2.54
$$
for any interval $J\subset [0,\,T].$
Then, using (2.54), we derive for any fixed \,
$0\leq \ell<\ell'\leq \lfloor T/\delta_\nu \rfloor$ \,
the bound
$$
\align
\Big|y_\nu^\ds(\ell'\delta_\nu)-y_\nu^\ds(\ell\delta_\nu)\Big|&=
\Big|y_\nu(\ell'\delta_\nu)-y_\nu(\ell\delta_\nu)\Big|
\\
\noalign{\medskip}
&\leq 
\tv\big\{y_\nu~;~[\ell\delta_\nu,\,\ell'\delta_\nu]\big\}
\\
\noalign{\medskip}
&\leq (1+C)\cdot(\ell'-\ell)\delta_\nu.
\tag 2.55
\endalign
$$
Therefore $y_\nu^\ds(\cdot)$ is a uniformly
bounded sequence
of Lipschitz maps, having Lipschitz constant \,
$\text{Lip}(y_\nu^\ds)\leq (1+C).$ \, Hence, applying Ascoli-Arzel\`a
Theorem, we can find a subsequence, that we still denote
$y_\nu^\ds(\cdot),$ which converges to some function \,
$y:[0,\,T]\mapsto \Bbb R^n,$ \,
uniformly on $[0,\,T].$ On the other hand,
by construction and thanks to (2.54), for any fixed $0\leq t\leq T,$ with 
$\ell \delta_\nu \leq t < (\ell+1)\delta_\nu,$ \ there holds
$$
\align
\Big|y_\nu(t)-y_\nu^\ds(t)\Big|&\leq
\Big|y_\nu(t)-y_\nu\big(\ell \delta_\nu\big)\Big|+
\Big|y_\nu\big(\ell \delta_\nu\big)-y_\nu^\ds(t)\Big|
\\
\noalign{\medskip}
&\leq \Big|y_\nu(t)-y_\nu(\ell \delta_\nu)\Big|+
\Big|y_\nu\big((\ell+1) \delta_\nu\big)-y_\nu\big(\ell \delta_\nu\big)\Big|
\\
\noalign{\medskip}
&\leq 2 \cdot
\tv\big\{y_\nu~;~[\ell \delta_\nu,\,(\ell+1)\delta_\nu]\big\}
\\
\noalign{\medskip}
&\leq 2 (1+C) \cdot \delta_\nu.
\tag 2.56
\endalign
$$
Thus, since $\delta_\nu \to 0$ as $\nu \to \infty,$
the uniform convergence 
of \, $y_\nu^\ds(\cdot)$ \, to \, $y(\cdot)$ \,
implies 
$$
\lim_{\nu\to\infty} 
\big\|y_\nu-y\big\|_{\L^\infty([0,T])} =0.
\tag 2.57
$$
By assumption, $\text{Range}(y_\nu)\subset A\subset \Omega$ \,
for every $\nu,$ and hence
from (2.57) we deduce that also the limit function $y(\cdot)$ takes values
inside~$\Omega$.  We can thus apply Theorem~1 to the
sequence \, $y_\nu(\cdot)$ \, and conclude that
the function $y : [0,\,T]\mapsto \Omega$ is a Carath\'eodory
solution of the unperturbed equation~(1.5)
with
$$
\big\|y-y_\nu\big\|_{\L^\infty([0,T])}
< \ve\,
$$
for all $\nu$ sufficiently large. We thus obtain a
contradiction with (2.52), concluding the proof.
\fine

\newpage

\vsk
\n{\medbf 3 - Robustness of Patchy Feedbacks.}
\v
In this section we apply the previous results on patchy vector fields
with impulsive perturbations, and 
construct (discontinuous) stabilizing feedback controls 
that enjoy robustness properties in the presence of 
measurement
errors and external disturbances.
Consider the nonlinear control system on~$\R^n$
$$
\dot x=f(x,u) \qquad\qquad u(t)\in \K,
\tag 3.1
$$
assuming that the control set $\K\subset\R^m$ is compact and that the map
$f:\R^n\times\R^m\mapsto\R^n$ is smooth.
We seek a 
feedback control $u=U(x)\in \K$ 
that stabilizes the trajectories of the closed-loop system
$$
\dot x=f\big(x,~U(x)\big)
\tag 3.2
$$
at the origin. It is well known that, 
even if every initial state \, $\overline x\in \Bbb R^n$ \,
can be steered to the origin by an open-loop control \, 
$u=u^{\overline x}(t),$ \, 
a topological obstruction can prevent the
existence of a continuous feedback control \, $u=U(x)$ \, which (locally)
stabilizes the system~(3.1).
This fact was first pointed out by Sussmann 
\Su \ for a two-dimensional system
($n=2,$ \ $\K= \Bbb R^2$), \ and by Sontag and Sussmann \SS \
for one-dimensional systems 
($n=1, \ \K=\Bbb R$). For general nonlinear systems, it
was further analyzed
by Brockett \Bro \ and Coron \
\Cornc.
It is thus natural to look for a stabilizing control  within
a class of discontinuous functions. However, this leads 
to a theoretical
difficulty, because, when the function \, $U$ \, 
is discontinuous, the differential
equation (3.2) may not admit any Carath\'eodory solution.
To cope with this problem, two different approaches have been
pursued.
\v
\n{\bf 1.} 
An algorithm is defined, which constructs approximate trajectories
in connection with an arbitrary (discontinuous)
feedback control function. For example, one can sample the feedback
control at a discrete set of times.
In this case, one is not concerned with
the existence of exact solutions, but only in the
asymptotic stabilization properties of all approximate solutions.

\v
\n {\bf 2.} Alternatively, by the asymptotic controllability
to the origin of system (3.1) by means of open-loop controls,
one proves the existence of a stabilizing feedback \, $u=U(x)$ \, 
having only a particular type of discontinuities.  
This feedback thus generates a {\it patchy vector field},
and the corresponding
system (3.2) always admits Carath\'eodory solutions.
\v

The first approach was initiated in \CLSS, 
and further developed in
\Riscl, \Risclsf.
The second was introduced in \AB,
defining the following class of piecewise constant
feedback controls:
\v
\n{\bf Definition 3.1.}  \ Let $\big(\Omega,\ g,\ 
(\Omega_\alpha,\,g_\alpha)_{_{\alpha\in \A}} \big)$ be a patchy
vector field. Assume that there exist control values
$k_\alpha\in \K$ such that, for each $\alpha\in\A,$ there holds
$$
g_\alpha(x) \doteq f(x,\, k_\alpha)\qquad\qquad\forall~  x \in 
D_\alpha \doteq 
\Omega_\alpha \setminus \bigcup_{\beta > \alpha} \Omega_\beta.
\tag 3.3
$$
Then, the piecewise constant map
$$
U(x) \doteq  k_\alpha\qquad \hbox{if}\qquad x \in D_\alpha
\tag 3.4
$$
is called a \ {\it patchy feedback} \ control on $\Omega,$
and referred to as \, 
$\big(\Omega,\ U,\ (\Omega_\alpha,\,k_\alpha)_{_{\alpha\in \A}} \big)$ \.
\v
\n{\bf Remark 3.2.}  \ From Definitions~1.2 and 3.1,
it is clear that
the field
$$
g(x)=f\big(x,\,U(x)\big)
$$
defined in connection with a given patchy feedback \
$\big(\Omega,\ U,\ (\Omega_\alpha,\,k_\alpha)_{_{\alpha\in \A}} \big)$ \
is precisely the patchy vector field \
$\big(\Omega,\ g,\ (\Omega_\alpha,\,g_\alpha)_{_{\alpha\in \A}} \big)$ \
associated with a family of fields $\big\{g_\alpha : \alpha\in \A\big\}$ 
satisfying (1.1) 
Clearly, the patches $(\Omega_\alpha,\,g_\alpha)$
are not uniquely determined by the patchy feedback $U$. Indeed, 
whenever $\alpha<\beta$, by (3.3) the values of $g_\alpha$ on the 
set~$\Omega_\alpha\setminus\Omega_\beta$ are irrelevant.
Moreover, recalling the notation (1.3) we have
$$
U(x) = k_{\alpha^\ast(x)}\qquad\quad \forall~x \in \Omega.
\tag 3.5
$$
\v

Here, we address the issue of robustness of a stabilizing
feedback law \, $u=U(x)$ \, w.~r.~t.~small 
internal and external perturbations
$$
\dot x=f\big(x,~U(x+\zeta(t))\big)+ d(t),
\tag 3.6
$$
where $\zeta=\zeta(t)$ represents a
state measurement error, and $d=d(t)$ represents
an external disturbance of the system dynamics (3.2).
Since we are dealing with a discontinuous
O.D.E., one cannot expect the full robustness
of the feedback $U(x)$ with respect to
measurement errors because of possible chattering
behaviour that may arise at discontinuity points
(see \Hdvf, \Ss).
Therefore, we shall consider state measurement
errors which are small in BV norm, avoiding such phenomena.
\v
Before stating our main result in this direction,
we recall here some basic definitions
and Proposition~4.2  in \AB\,. This provides
the semi-global practical stabilization 
(steering all states from a given compact set of initial data 
into a prescribed neighborhood of zero) of
an asymptotycally controllable system, by means of a
patchy feedback control which is robust with respect
to external disturbances.
We consider as (open-loop) 
{\it admissible controls} all the measurable functions
$u : [0,\,\infty)\to \R^m$ such that $u(t) \in \K$  for a.e. 
$t\geq 0$. \ 
\v
\n{\bf Definition 3.3.}  \ The system (3.1) is
globally {\it asymptotically controllable} to the origin 
if the following 
holds.
\vskip 0.5em
\item{\bf 1.}  \ {\bf Attractiveness:} \ for each 
$\overline x\in \R^n$ there exists
some admissible (open-loop)  
control $u =  u^{\overline x}(t)$ such that the 
corresponding trajectory of
$$
\dot x(t) = f\big(x(t), \, u^{\overline x}(t)\big),
\qquad\quad x(0)=\overline x\,,
\tag 3.7
$$
either reaches the origin in finite time,
or tends to the origin as $t \to \infty.$
\v
\item{\bf 2.}  \ {\bf Lyapunov stability:} \ for each $\eps>0$ there 
exists $\delta>0$
such that the following holds. 
For every $\overline x\in\R^n$ with $|\overline x|<\delta,$ 
there is an admissible
control $u^{\overline x}$ as in 1. steering the system
from $\overline x$ to the origin,
so that the corresponding trajectory
of (3.7) satisfies $|x(t)|<\eps$ for all $t\geq 0.$
\v
\n{\bf Proposition 3.4.} [A-B, \, Proposition 4.1] \
{\it Let system (3.1) be globally asymptotically
controllable to the origin. Then, 
for every \ $0<r<s,$ \ 
one can find \
$T>0,\  
\chi>0,$ \ and a patchy feedback 
control \, $U : D\mapsto \K,$ \,
defined on some domain 
%
$$
D \supset 
\big\{x\in \Bbb R^n~;~ r \leq |x|\leq s \big\}
%
\tag 3.8
$$
so that the following holds. 
For any measurable map \, $d : [0,\, T]
\mapsto \R^n$ \, such that
$$
\big\|d\big\|_{\L^\infty([0,\,T])}\leq \chi\,,
$$ 
and for any initial state \, $x_0$ \,
\, with \, $r \leq |x_0| \leq s,$ \,
the perturbed system 
$$
\dot x = f\big(x,\,U(x)\big) + d(t)
\tag 3.9
$$
admits a (forward) Carath\'eodory trajectory starting from~$x_0.$
Moreover, for any such trajectory \,
$t \mapsto \gamma(t),$ \, $t\geq 0,$ \,
one has
$$
\gamma(t)\in D
\qquad \quad \forall~t \geq 0,
\tag 3.10
$$
and 
there exists \,
$\overline t_{\gamma}<T$ \,
such that
$$
\big|\gamma(\,\overline t_{\gamma})\big| < r.
\tag 3.11
$$
}
\v

Relying on Corollary~1.3 of Theorem~1 and on
Proposition 3.4., we obtain
here the following result concerning robustness
of a stabilizing feedback w.~r.~t.~both
internal and external perturbations.
\vs

\n{\bf Theorem 2.} {\it 
Let system (3.1) be globally asymptotically
controllable to the origin. 
Then, for every
\ $0<r<s,$ \ one can find \
$T'>0,\ 
\chi'>0,$ \ 
and a patchy feedback 
control \,
$U' : D'\mapsto \K$ \,
defined on some domain \, $D'$ \, 
satisfying (3.8), so that the following
holds.
Given any pair of maps \, $\zeta \in BV([0,\,T']),$ 
$d \in \L^\infty([0,\,T']),$ \, such that
$$
\tv\{\zeta\}\leq \chi'\,,
\qquad\quad
\big\|d\big\|_{\L^\infty([0,\,T'])}\leq \chi'\,,
\tag 3.12
$$ 
and any initial state \, $x_0$ \, with \, $r \leq |x_0| \leq s,$ \, 
for every solution \,
$t \mapsto x(t),$ \, $t\geq 0,$ \,
of the perturbed system~(3.6) starting from $x_0,$ 
one has
$$
x(t)\in D'
\qquad\ \forall~t\in [0,\,T']\,,
\tag 3.13
$$
and there exists \
$\overline t_{x}<T',$ \,
such that
$$
\big|x(\,\overline t_{x})\big| < r.
\tag 3.14
$$
}
\v

\noindent
{\bf Proof.}  \

\noindent
{\smbf 1.} \ 
Fix $0<r<s.$ \, Then, according with Proposition~3.4,
we can find
\ $T'>0$,\ 
and a patchy feedback 
control \, $U' : D'\mapsto \K,$ \,
defined on some domain 
%
$$
D'\supset 
\big\{x\in \Bbb R^n~;~ r/3 \leq |x|\leq s \big\}\,,
\tag 3.15
$$
%
so that the following holds. 
For every Carath\'eodory solution \,
$t \mapsto x(t),$ \, $t\geq 0$ \,
of the unperturbed system~(3.2) (with \, $U=U'$ \,) starting from  a point
\, $x_0$ \, in the compact set \, 
$$
K\doteq
\big\{x\in \Bbb R^n~;~ r \leq |x|\leq s \big\},
\tag 3.16
$$
one has
$$
x(t)\in D_\rho\doteq
\big\{x\in D'~:~ d\big(x,\, \partial D'\big)> \rho\big\}
\qquad \quad \forall~t \geq 0,
\tag 3.17
$$
for some constant \, $\rho>0.$ \,
Moreover, there exists \
$\overline t_{x}<T'$ \,
such that
$$
\big|x(\,\overline t_{x})\big| < \frac{r}{3}.
\tag 3.18
$$
According with Definition~3.1, the field
$$
g(x)\doteq f\big(x, \, U'(x)\big)
\tag 3.19
$$
is a patchy vector field associated to the family
of fields \, $\big\{g_\alpha : \alpha\in \A\big\}$ \,
defined as in (3.3). 
The smoothness of \, $f$ \, guarantees that,
for BV perturbations \, $w=w(t)$ \, having 
some uniform bound \, $\tv \{w \}\leq \widehat \chi$ \,
on the total variation,
every (left-continuous) solution  \,
 $y : [0,T']\mapsto\Bbb R^2$ \, of the impulsive equation~(1.7),
starting at a point \, $x_0\in K,$ \, takes values in the
closed set
$$
A\doteq B(D_\rho,\, \rho/2).
\tag 3.20
$$
Therefore, thanks to Corollary~1.3 of Theorem~1, 
there exists some constant
$$
0<\widehat\chi'=\widehat\chi'(A,\, K,\, T',\, r/3)<\widehat\chi
\tag 3.21
$$
such that the following holds.
If $y:[0,T']\mapsto \Bbb R^2$ is a (left-continuous)
solution of the impulsive equation (1.7),
with $y(0)\in K$ and $\tv (w)<\widehat\chi'$, 
then one has 
$$
y(t)\in A
\qquad \quad \forall~t\in [0,\,T']\,,
\tag 3.22
$$
and there exists \
$\overline t_{y}<T'$ \,
such that
$$
\big|y(\,\overline t_{y})\big| < \frac{2 r}{3}.
\tag 3.23
$$
\vs

\noindent
{\smbf 2.} \ In connection with the patchy feedback \, $U'$ \, introduced
above, define the map
$$
h(y,z) 
\doteq f\big(y-z,\, U'(y)\big)- 
f\big(y,\, U'(y)\big) 
\tag 3.24
$$
and observe that, by the smoothness of \, $f,$ \,
there will be some constant \, $\overline c>0,$ \,
such that 
$$
\big|h(y,z)\big| \leq \overline c\cdot |z|
\qquad\quad \forall~y\in A,\quad |z|\leq \widehat \chi'.
\tag 3.25
$$
Consider now a 
pair of maps \, $\zeta \in BV([0,\,T']),$ 
$d \in \L^\infty([0,\,T']),$ \, satisfying (3.12)
with
$$
\chi'<\min\left\{
\frac{\widehat \chi'}{2 (1 + T' \overline c')},\ \frac{r}{3}
, \ \frac{\rho}{2}
\right\}\,,
\tag 3.26
$$
and  
let \, $x=x(t)$ \, be any Carath\'eodory
solution of the perturbed system (3.6), 
with an initial condition \, $x(0)=x_0\in K.$ \,
Then, as observed in the introduction, the map \,
$$
t \mapsto y(t)\doteq x(t)+\zeta(t)
\tag 3.27
$$ 
satisfies the impulsive equation (1.7) where
$$ 
w (t)\doteq \zeta(t)+\int_0^t 
\big(h\big(y(s), \zeta(s)\big)  +d(s)\big)\,ds\,.
\tag 3.28
$$
But then, since (3.12), (3.25), (3.26), together, imply
$$
\align
\tv\big\{ w~;~ [0,\, T']\big\}
&\leq \tv\big\{ \zeta~;~ [0,\, T']\big\}+ T'
\overline c \cdot \big\|\zeta\big\|_{\L^\infty([0,\,T'])}+
T'\cdot \big\|d\big\|_{\L^\infty([0,\,T'])}
\\
\noalign{\medskip}
&\leq \big(1 + T'\overline c\big)\cdot \tv\big\{ \zeta~;~ [0,\, T']\big\}+
\big\|d\big\|_{\L^\infty([0,\,T'])}
\\
\noalign{\medskip}
&<\widehat\chi'\,,
\endalign
$$
from (3.22)-(3.23) and (3.12),  (3.20), (3.26), (3.27) it follows
$$
\gathered
x(t)\in B(A,\, \rho/2)\subset D'\,,
\qquad\ \forall~t\in [0,\,T']\,,
\\
\noalign{\medskip}
\big|x(\,\overline t_{y})\big|< r
\qquad\ \text{for \ some}\qquad \overline t_{y}<T'\,,
\endgathered
\tag 3.29
$$
which completes the proof of the theorem, taking \,
$\chi'$ \, as in (3.26).
\fine
\vs

\v
\n{\bf Remark~3.5.}  \  
For discontinuous stabilizing feedbacks constructed 
in terms of sampling solutions, 
an alternative concept of robustness 
was introduced in \CLRSfslf, \Ss.
In this case, one considers a
partition of the time interval 
and 
applies a constant control
between two consecutive sampling times.
To preserve stability,
the measurement error should be 
sufficiently small compared to the maximum step size.
Moreover, each step size 
should be big enough
to prevent possible
chattering phenomena.
The next result shows that the feedback provided by
[A-B, \, Proposition 4.2]
enjoys also this type of robustness.
Before stating this result we describe now the concept
of sampling trajectory associated
to the perturbed system (3.6)
that was introduced in \CLRSfslf.
\vskip 0.8pt

Let  an initial condition \, $x_0,$ and a partition
 \, $\pi=\{0=\tau_0< \tau_1< \cdots < \tau_{m+1}=T\}$ \,
of the interval~$[0,\, T],$ \, be given
A sampling trajectory \, $x_\pi$ \,
of the perturbed system (3.6),
corresponding to
a set of measurement errors \, $\{e_i\}_{i=0}^m$ \,
and an external disturbance \, $d \in \L^\infty([0,\,T]),$ \,
is defined in a step-by-step
fashion as follows.
Between \, $\tau_0$ \, and \, $\tau_1,$ \, let $x_\pi(\cdot)$ \,
be a Carath\'eodory solution of
$$
\dot x=f\big(x,~U(x_0+e_0)\big)+ d(t)
\qquad\ t\in [\tau_0,\, \tau_1]\,,
\tag 3.30
$$
with initial condition \, $x_\pi(0)=x_0.$ \,
Then, \, $x_\pi(\cdot)$ \,
is recursively obtained by solving the system
$$
\dot x=f\big(x,~U(x_\pi(\tau_i)+e_i)\big)+ d(t)
\qquad\ t\in [\tau_i,\, \tau_{i+1}],\quad i>0.
\tag 3.31
$$
The sequence \, $\{x_\pi(\tau_i)+e_i\}_{i=0}^m$ \,
corresponds to the non-exact measurements
used to select control values.
\vs

\n{\bf Theorem 3} \ {\it 
Let system (3.1) be globally asymptotically
controllable to the origin. 
Then, for every
\ $0<r<s,$ \ one can find \
$T''>0,$ \, 
$\chi''>0,$ \, 
$\overline\delta>0,$ \,
\, $\overline k>0,$ \,
and a patchy feedback 
control \,
$U'' : D''\mapsto \K$ \,
defined on some domain \, $D''$ \,
satisfying (3.8)
so that the following holds.
Given an initial state \, $x_0$ \, with \, $r \leq |x_0| \leq s,$ \,
a partition 
 \, $\pi=\{\tau_0=0,\, \tau_1,\, \dots , \tau_{m+1}=T''\}$ \,
of the interval \, $[0,\, T''\,]$ \,
having the property
$$
\hskip 1in
\frac{\delta}{2}\leq \tau_{i+1}-\tau_i \leq \delta \qquad\forall~i\,,
\qquad\quad
\text{for \ some}\qquad
\delta \in ~]0,\, \overline \delta\,]\,,
\tag 3.32
$$
a set of measurement errors \, $\{e_i\}_{i=0}^m\,$ \,
and an external disturbance \, $d \in \L^\infty([0,\,T''\,])$ \,
that satisfy
$$
\align
\max_{i} |e_i| 
&\leq \overline k \cdot \delta\,, 
\tag 3.33
\\
\noalign{\medskip}
\big\|d\big\|_{\L^\infty}
&\leq \chi''\,,
\tag 3.34
\endalign
$$
the resulting sampling solution \,
$x_\pi(\cdot)$ \, 
starting from $x_0$ 
has the property
$$
x_\pi(t)\in D''
\qquad\ \forall~t\in [0,\,T''\,]\,.
\tag 3.35
$$
Moreover, there exists \
$\overline t_{x_\pi}<T'',$ \,
such that
$$
\big|x_\pi(\,\overline t_{x_\pi})\big| < r.
\tag 3.36
$$
}
\v

\noindent
{\bf Proof.}  \

\noindent
{\smbf 1.} \ Fix \, $0<r<s.$ \,
Then,  according with Proposition~3.4,
we can find \,
$
T'>0, \,
\chi'>0, \,
$
and a patchy feedback control \,
$
U'' : D''
\longmapsto \K
$ \,
defined on a domain 
%
$$
D''\supset
\big\{x\in \Bbb R^n~;~ r/3 \leq |x|\leq 2s \big\}
$$
so that the following holds.
For every external disturbance \, $d \in \L^\infty$  \,
satisfying (3.34) with \, $\chi'' \leq   \chi'$\,,
%
and for any Carath\'eodory solution \,
$t \mapsto x(t),$ \, $t\geq 0$ \,
of the perturbed system~(3.9) (with \, $U=U''$), starting from  a point
\, $x_0$ \, with \, $r \leq |x_0| \leq s,$ \,
one has
$$
x(t)\in D_{\rho_1}\doteq
\big\{x\in D''~:~ d\big(x,\, \partial D''\big)> \rho_1\big\}
\qquad \quad \forall~t \geq 0,
\tag 3.37
$$
for some constant \, $\rho_1>0.$ \
Moreover, there exists \
$\overline t_{x}<T'$ \,
such that
$$
\big|x(\,\overline t_{x})\big| < \frac{r}{3}.
\tag 3.38
$$
Let 
$$
\big\{(\Omega_\alpha,\ g_\alpha)~:~\alpha=1,\dots,N\big\}
\qquad\quad g_\alpha(x)=f(x,k_\alpha),\qquad k_\alpha\in \K,
\tag 3.39
$$
be the collection of patches associated with the patchy vector field
$$
g(x)=f\big(x,\, U''(x)\big).
\tag 3.40
$$
We may assume  
that every vector field \,$ g_\alpha$ \, 
is
defined on a  neighborhood 
\, $B(\Om_\alpha,\, \rho_2), \ 0<\rho_2 \leq \rho_1$ \,
of the domain \, $\Om_\alpha$ \, 
so that, setting
$$
\Om_\alpha^\rho \doteq 
\big\{x \in \Om_\alpha~;~ d\big(x,\ \partial \Om_\alpha\big)>\rho\big\}\,,
\tag 3.41
$$
one has
$$
\Om_\alpha^{\rho_2}\neq \emptyset\,,
$$
and that every \, $g_\alpha$ \,
is uniformly non-zero on the domain \, $D_\alpha$ \, 
defined in (3.3).
Moreover, thanks to the inward-pointing condition (1.1),
we may choose the constants \, $\ 0<\rho_2 < r/3,$ \, 
and \, $\chi''\leq \chi'$ \, so that there holds
$$
\big|g_\alpha(x)\big|\geq 2 \chi''
\qquad\quad \forall~ x\in B(D_\alpha,\, \rho_2),
\tag 3.42
$$
and
$$
\la g_\alpha(x) + v,~\bn(x)\ra <0\qquad\forall~x
\in B(\partial\Om_\alpha,\, \rho_2),
\qquad |v|\leq \chi''\,.
\tag 3.43
$$
%
For every \, $d \in \L^\infty,$ \,
we denote
by \, $t \mapsto x^\alpha\big(t;\, t_0,\, x_0,\ d\big)$ \,
the solution of the Cauchy problem
%
$$
\dot x = g_\alpha(x)+d(t),
\qquad
x(t_0)=x_0,
\tag 3.44
$$
and let \, $[t_0, \, t^{\max}]$ \, be the domain of definition 
of the maximal (forward)
solution of (3.45) that is
contained in \, $B(D_\alpha,\, \rho_2).$ \,
\v

Observe that, 
since every Carath\'eodory solution of the 
perturbed system~(3.9) (with $U=U''$), starting from a point 
\, $x_0\in
B(0,s)\!\setminus\!\overset \,\circ \to B(0,\,r),$ \
reaches the interior of the ball \, $B(0,\, r/3)$ \, in finite time,
and because of (3.42),
for any \,$\alpha = 1, \dots ,N$ \, one can find \, $T_\alpha>0$ \,
with the following property.
\vvs
\parindent 30pt
\item{\bf (P)$_{\bold 1}$} \
For every 
\, $x_0\in B(D_\alpha, \, \rho/2),$ \, \, $0<\rho< \rho_2,$ \, 
and for any \, $d \in \L^\infty$ \,
satisfying (3.34),
there exists some time \, $t_\rho\doteq t_\rho(x_0,\, d)< T_\alpha$ \,
such that, either one has
$$
\big|x^\alpha\big(t_0+t_\rho;\, t_0,\, x_0,\, d\big)\big|<\frac{2r}{3},\,
\tag 3.45
$$
or else there holds 
$$
x^\alpha\big(t;\, t_0,\, x_0,\, d\big)\in 
B(D_\alpha,\, \rho_2)\setminus B(D_\alpha, \, \rho)
\qquad\quad \forall~t\in[t_0+t_\rho,\ t^{\max}].
\tag 3.46
$$
\parindent 20pt

\vvs
\noindent
On the other hand, relying on the inward-pointing condition (3.43), 
we deduce two further properties of the solutions of (3.44).
\vvs
\parindent 30pt
\item{\bf (P)$_{\bold 2}$} \
The sets \, $\Om_\alpha^{\rho},$ \, $0<\rho\leq  \rho_2,$ \, defined in (3.41)
are positive
invariant regions for trajectories of~(3.44), i.e.,
for every 
\, $x_0\in \Om_\alpha^\rho,$ \, 
and for any \, $d \in \L^\infty$ \,
satisfying (3.34),
one has
$$
x^\alpha\big(t;\, t_0,\, x_0,\, d\big)\in 
\Om_\alpha^{\rho}
\qquad\quad \forall~t \geq t_0\,.
\tag 3.47
$$
\parindent 20pt

\vvs
\parindent 30pt
\item{\bf (P)$_{\bold 3}$} \
There exists some constant \, $\overline c>0$ \, so that,
for every 
\, $x_0\in B(\Om_\alpha, \, \rho),$ \, \, $0<\rho \leq  \rho_2,$ \, 
such that \, $d\big(x_0,\ \partial \Om_\alpha\big) \leq \rho,$ \,
and for any \, $d \in \L^\infty$ \,
satisfying (3.34),
one has
$$
x^\alpha\big(t;\, t_0,\, x_0,\, d\big)\in 
\Om_\alpha^{2 \rho}
\qquad\quad \forall~t \geq t_0+\overline c \cdot \rho\,.
\tag 3.48
$$
\parindent 20pt
\vs

\noindent
{\smbf 2.} \ Consider an initial state \, $x_0\in 
B(0,s)\!\setminus\!\overset \,\circ \to B\big(0,\,r\big),$ \, 
and a partition 
 \, $\pi=\{\tau_i\}_{i\geq 0}$ \, of \, $[0,\, \infty[$ \,
having the property (3.32), with 
$$
\gather
0<\delta \leq 
\overline \delta \doteq \min\Big\{\overline c\cdot \rho_2, \ 
\frac{\rho_1}{M}\Big\},
\tag 3.49
\\
\noalign{\medskip}
M\doteq 
\sup\big\{|g_\alpha(x)|~:~x\in B(\Om_\alpha,\, \rho_2),
\quad \alpha=1,\dots,N\big\}\,.
\endgather
$$
%
%
Let $x_\pi: [0,\, \infty[ ~\mapsto \Bbb R^n$ \,be a sampling solution 
starting from \, $x_0,$ \, and corresponding to 
a set of measurement errors \, $\{e_i\}_{i=0}^m\,$ \,
and to an external disturbance \, $d(\cdot) \in \L^\infty$ \,
that satisfy (3.33)-(3.34)
with 
$$
\overline k \doteq \frac{1}{2 \overline c}\,.
\tag 3.50
$$
We will first show the following
\v

\noindent
{\bf Lemma 3.6.} {\it 
The map
$$
i \longmapsto \alpha^*(\tau_i) \doteq \alpha^*\big(x_\pi(\tau_i)+e_i\big)
\qquad\ i \geq 0\,,
\tag 3.51
$$
is non-decreasing.} 
\v

\noindent
Indeed, assume that  \, $\alpha^*(\tau_i)=\widehat \alpha,$ \, 
which, by definitions (1.3), (3.3), (3.5)
implies  
$$
\gather
x_\pi(\tau_i)+e_i \in D_{\widehat \alpha}\,,
\tag 3.52
\\
\noalign{\medskip}
x_\pi(\tau_{i+1})= x^{\widehat \alpha}\big(\tau_{i+1};\ 
\tau_i,\, x_\pi(\tau_i),\, d\restriction_{[\tau_i,\, \tau_{i+1}]}\big)\,,
\tag 3.53
\endgather
$$ 
Then,
because of (3.33), (3.49)-(3.50),
one has
$$
x_i \doteq x_\pi(\tau_i)
\in B(D_{\widehat\alpha},\ \overline k \delta) \subset 
B(\Om_{\widehat\alpha},\ \rho_2).
\tag 3.54
$$
We shall consider separately the case in which
$$
x_i\in D_{\widehat \alpha}^{\overline k \delta}
\subset \Om_{\widehat \alpha}^{\overline k \delta}
\qquad\quad \overline k \delta\leq \rho_2\,,
\tag 3.55
$$
and the case where
$$
x_i\in B(D_{\widehat\alpha},\ \overline k \delta),
\qquad\quad d\big(x_i\,,\, \partial \Om_{\widehat\alpha}\big)\leq 
\overline k \delta\leq \rho_2\,.
\tag 3.56
$$
In the first case, 
using (3.53) and 
applying {\bf (P)$_{\bold 2}$} we deduce that \, 
$x_\pi(\tau_{i+1})\in \Om_{\widehat \alpha}^{\overline k \delta}$ \,
which, in turn, because of (3.33), (3.49)-(3.50), implies
$$
x_\pi(\tau_{i+1})+e_{i+1}\in \Om_{\widehat \alpha}.
\tag 3.57
$$
{}From (3.57), by definition (1.3) we derive \, 
$$
\alpha^*(\tau_{i+1})\geq \widehat \alpha\,,
\tag 3.58
$$
proving the lemma whenever (3.55) holds. On the other hand,
when (3.56) is verified, since
by (3.32), (3.50) one has \, 
$$
\tau_{i+1}-\tau_i \geq \frac{\delta}{2} = 
\overline c \overline k \cdot \delta\,,
$$
applying
{\bf (P)$_{\bold 3}$} we deduce \, 
$x_\pi(\tau_{i+1})\in \Om_{\widehat \alpha}^{2\overline k \delta}$.
This again implies (3.57)-(3.58), completing the proof of Lemma 5.6.
\v

Next, relying on {\bf (P)$_{\bold 1}$}, and
setting
$$
\hskip 0.7in
\aligned
i'_\alpha
&\doteq \min \big\{\, i\geq 0 \ \ ;\ \ \alpha^*({\tau_i}) = 
\alpha\,, \qquad 
x_\pi(\tau_i)\notin B(0,\, 2r/3)\,\big\}\,,
\\
\noalign{\medskip}
i''_\alpha
&\doteq \max \big\{\, i\geq 0 \ \ ;\ \ \alpha^*({\tau_i}) = 
\alpha\,, \qquad 
x_\pi(\tau_i)\notin B(0,\, 2r/3)\,\big\}\,,
\endaligned
\qquad \alpha \in \text{Range}(\alpha^*),
\tag 3.59
$$
%
we deduce
$$
\tau_{_{i''_\alpha}}-\tau_{_{i'_\alpha}} \leq T_\alpha
\qquad\quad \forall
~\alpha \in \text{Range}(\alpha^*)\,.
\tag 3.60
$$
Indeed, if (3.60) does not hold,  
by definitions (3.3), (3.5) one has
$$
\gather
x_{i'_\alpha}\doteq x_\pi(\tau_{_{i'_\alpha}})
\in B(D_\alpha,\ \overline k \delta) \subset 
B(\Om_\alpha,\ \rho_2/2)\,,
\tag 3.61
\\
\noalign{\medskip}
x_\pi(t)= x^\alpha\big(t;\ 
\tau_{_{i'_\alpha}},\, x_{i'_\alpha},\, 
d\restriction_{[\tau_{_{i'_\alpha}},\, \tau_{_{i''_\alpha+1}}]}\big)\,,
\qquad\ \forall~t\in [\tau_{_{i'_\alpha}},\ \tau_{_{i''_\alpha+1}}]\,.
\tag 3.62
\endgather
$$
But then,  applying {\bf (P)$_{\bold 1}$}, one could find some \,
$\widehat i \leq i''_\alpha$ \, such that
$$
x_\pi(t)\in 
B(D_\alpha,\, \rho_2)\setminus B(D_\alpha, \, 2 \overline k \delta)
\qquad\quad \forall~t\in [\tau_{_{\widehat i}}, \, \tau_{_{i''_\alpha+1}}]\,.
$$
By definitions (1.3), (3.51) and because of (3.33),
this implies
$$
\alpha^*(\tau_i)>\alpha \qquad \forall \widehat i \leq i \leq i''_\alpha\,,
$$
providing a contradiction with (3.59).
\v

To conclude the proof of Theorem 3, we observe that
the monotonicity of the map (3.51), together with the
estimate (3.60), implies that there exists some time 
\, $\overline t_{x_\pi}<T'' \doteq \sum_{\alpha=1}^N T_\alpha$ \,
such that (3.36) is verified.
Moreover, (3.35) clearly follows from 
(3.37) and (3.49).
\fine
\vs

\n{\bf Remark~3.7.}  \ Consider a partition 
 \, $\pi=\{\tau_0=0,\, \tau_1,\, \dots , \tau_{m+1}=T\}$ \,
of the interval \, $[0,\, T\,]$ \,
having the property (3.32).   If
we associate to a set of measurement errors \, $\{e_i\}_{i=1}^m$ \,
satisfying (3.33) the piecewise constant function
\, $\zeta~:[0,\,T]\mapsto \Bbb R^n$ \, defined as
$$
\zeta(t)=e_i \qquad \forall t\in ~]\tau_i,\, \tau_{i+1}]\,,
$$
then
$$
\tv \{\zeta\} \leq 4 \overline k \cdot T\,.
$$
Thus, taking the constant \, $\overline k$ \, sufficiently small
we may reinterpret the {\it discrete} internal disturbance 
allowed for a sampling solution in Theorem~3 as a particular case of 
the measurement errors with small total variation considered in Theorem~2.
\vs


%
\parindent 60pt

\vsk
\centerline{\medbf References.}
\vs
\noindent\item{\AB\ } {\smc F. Ancona, A. Bressan},
Patchy vector fields and asymptotic stabilization,
{\it ESAIM - Control, Optimiz. Calc. Var.,} 
Vol. {\bf 4}, (1999), pp. 445-471.
\medskip
%
%
%
\noindent\item{\Bimp\ } {\smc A. Bressan}, 
On differential systems with impulsive controls, 
{\it Rend. Sem. Mat. Univ. Padova}, 
Vol. {\bf 78}, (1987), pp.  227-235.
\medskip
\noindent\item{\Bro\ } {\smc  R.W. Brockett,}
Asymptotic stability and feedback stabilization, in \
{\it Differential \ Geometric \ Control Theory} \ (R.W. Brockett,
R.S. Millman, and H.J. Sussmann, eds.), \ Birkhauser, Boston, 
(1983), pp. 181-191.
\medskip
\noindent\item{\CLRSfslf\ } {\smc F.H. Clarke, Yu.S. Ledyaev, L. Rifford,
R.J. Stern,}
Feedback stabilization and Lyapunov functions,
{\it SIAM J. Control Optim.},
{\bf 39}, (2000), no. 1, pp. 25-48.
\medskip
\noindent\item{\CLSS\ } {\smc F.H. Clarke, Yu.S. Ledyaev, E.D. Sontag,
A.I. Subbotin,}
Asymptotic controllability implies feedback stabilization,
{\it IEEE Trans. Autom. Control}, {\bf 42} (1997), pp. 1394-1407.
\medskip
%
  
%
%
\noindent\item{\Cornc\ } {\smc  J.-M. Coron,} 
A necessary condition for feedback stabilization, 
{\it Systems Control Lett.}, 
{\bf 14} (1990), pp. 227-232.
\medskip
%
%
%
\noindent\item{\Cortv\ } {\smc  J.-M. Coron,} 
Stabilization in finite time of locally controllable systems 
by means of continuous \ time-varying \ feedback laws, \
{\it SIAM J. Control Optim.}, 
{\bf 33} (1995), \ pp.~804-833.
\medskip
%
%
%
%
%
%
\noindent\item{\Hdvf\ } {\smc  H. Hermes,} 
Discontinuous vector fields and feedback control, in 
{\it Differential Equations and Dynamical Systems}, 
(J.K. Hale and J.P. La Salle eds.), Academic Press, 
New York, (1967), pp.~155-165.
\medskip
%
%
%
%
\noindent\item{\KS\ } {\smc  N.N. Krasovskii and A.I. Subbotin,} 
Positional Differential Games, Nauka, Moscow, (1974) [in Russian].
Revised English translation: Game-Theoretical Control Problems,
Springer-Verlag, New York, 1988.
\medskip
\noindent\item{\LSrrs\ } {\smc Yu.S. Ledyaev and E.D. Sontag}, 
A remark on robust stabilization of general asymptotically
controllable systems, 
in {\it Proc. Conf. on Information Sciences and Systems (CISS 97)}, 
Johns Hopkins, Baltimore, MD, (1997), pp. 246-251.
\medskip
\noindent\item{\LScrs\ } {\smc Yu.S. Ledyaev and E.D. Sontag}, 
A Lyapunov characterization of robust stabilization, 
{\it Journ. Nonlinear Anal.}, {\bf 37} (1999), 
pp. 813-840.
\medskip
%
%
%
\noindent\item{\Riscl\ } {\smc L. Rifford}, 
Existence of Lipschitz and semiconcave control-Lyapunov
functions, 
{\it SIAM J. Control Optim.}, {\bf 39}
(2000), no. 4, pp. 1043-1064.
\medskip
\noindent\item{\Risclsf\ } {\smc L. Rifford}, 
Semiconcave control-Lyapunov
functions and stabilizing feedbacks, (2000),
preprint.
\medskip
\noindent\item{\Ry\ } {\smc E.P. Ryan}, 
On Brockett's condition for smooth stabilizability
and its necessity in a context of nonsmooth feedback, 
{\it SIAM J. Control Optim.}, 
{\bf 32} (1994), \ pp.~1597-1604.
\medskip
%
%
%
%
\noindent\item{\Ss\ } {\smc E.D. Sontag}, 
Stability and stabilization: discontinuities and the
effect of disturbances, in \,
{\it Proc. NATO Advanced Study Institute - 
Nonlinear Analysis, Differential Equations, and Control}, 
(Montreal, Jul/Aug 1998), F.H. Clarke and R.J. Stern eds.,
Kluwer, (1999),\  pp. 551-598.
\medskip
\noindent\item{\SS\ } {\smc E.D. Sontag and H.J. Sussmann}, 
Remarks on continuous feedback,  in
{\it Proc. IEEE Conf. Decision and Control}, Aulbuquerque,
IEEE Publications, Piscataway,
(1980), \ pp.~916-921.
\medskip
%
%

\end